\documentclass[a4paper,12pt]{article}
\usepackage{times}
\usepackage{amsmath}
\usepackage{amssymb}
\usepackage{amsthm}
\usepackage{graphicx}
\begin{document}
\title{Deconvolving oscillatory transients\\with a Kalman filter}
\author{Andreas M\"uller\footnote{
University of Applied Sciences, Oberseestrasse 10, CH-8640 Rapperswil,
Switzerland, Email: {\tt andreas.mueller@hsr.ch}}}
\date{}
\maketitle
\begin{abstract}
This paper describes a method to filter oscillatory transients from
measurements of a time series which were at least an order of magnitude
larger than the signal to be measured. Based on a Kalman filter,
it has an optimality property and
a natural scaling parameter that allows to tune it
to high resolution or low noise.
\end{abstract}
\begin{figure}
\begin{center}
\includegraphics[width=0.8\hsize]{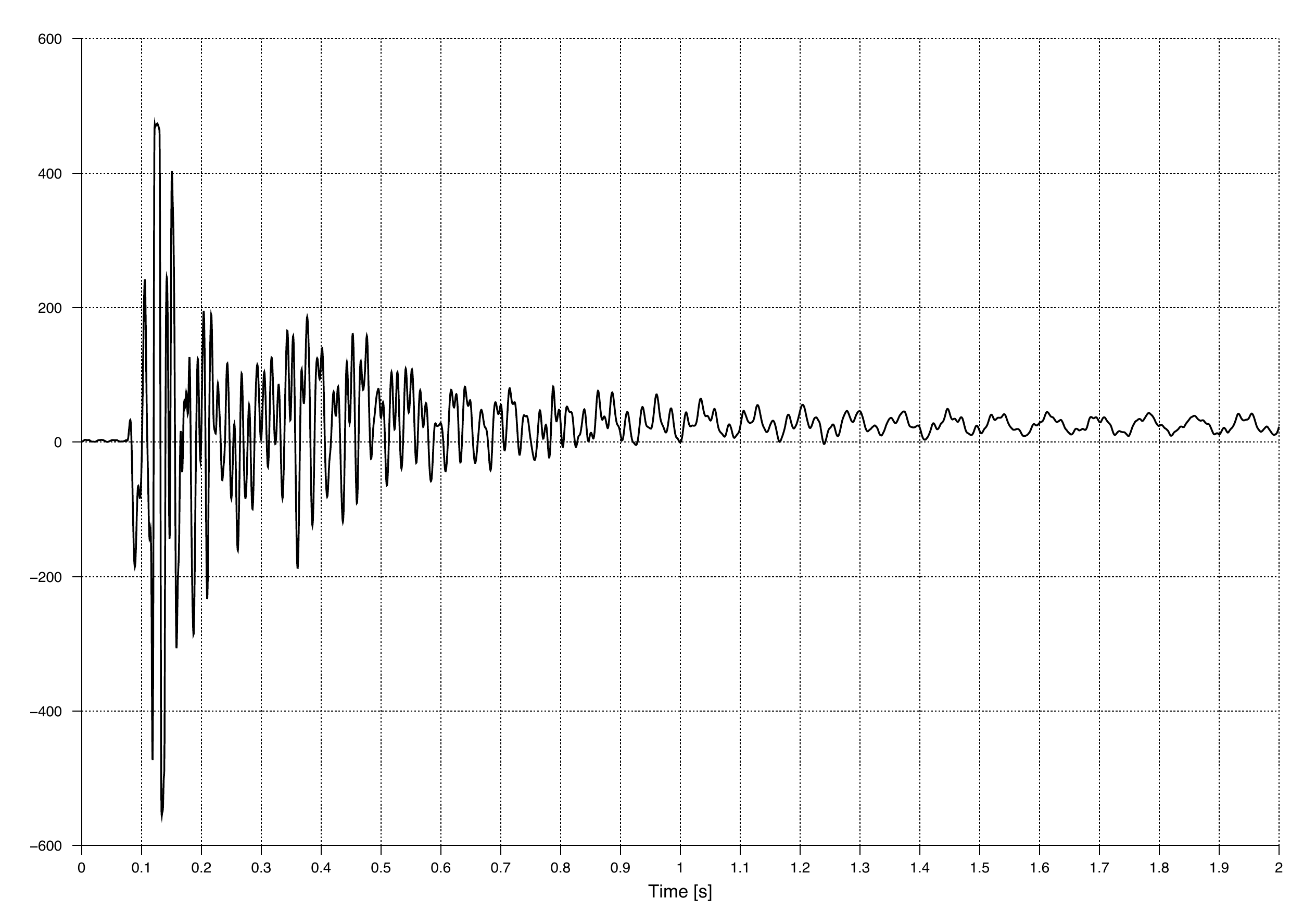}
\end{center}
\caption{Wind tunnel measurements of opening shock dominated by
transients\label{raw}}
\end{figure}
\begin{figure}
\begin{center}
\includegraphics[width=0.8\hsize]{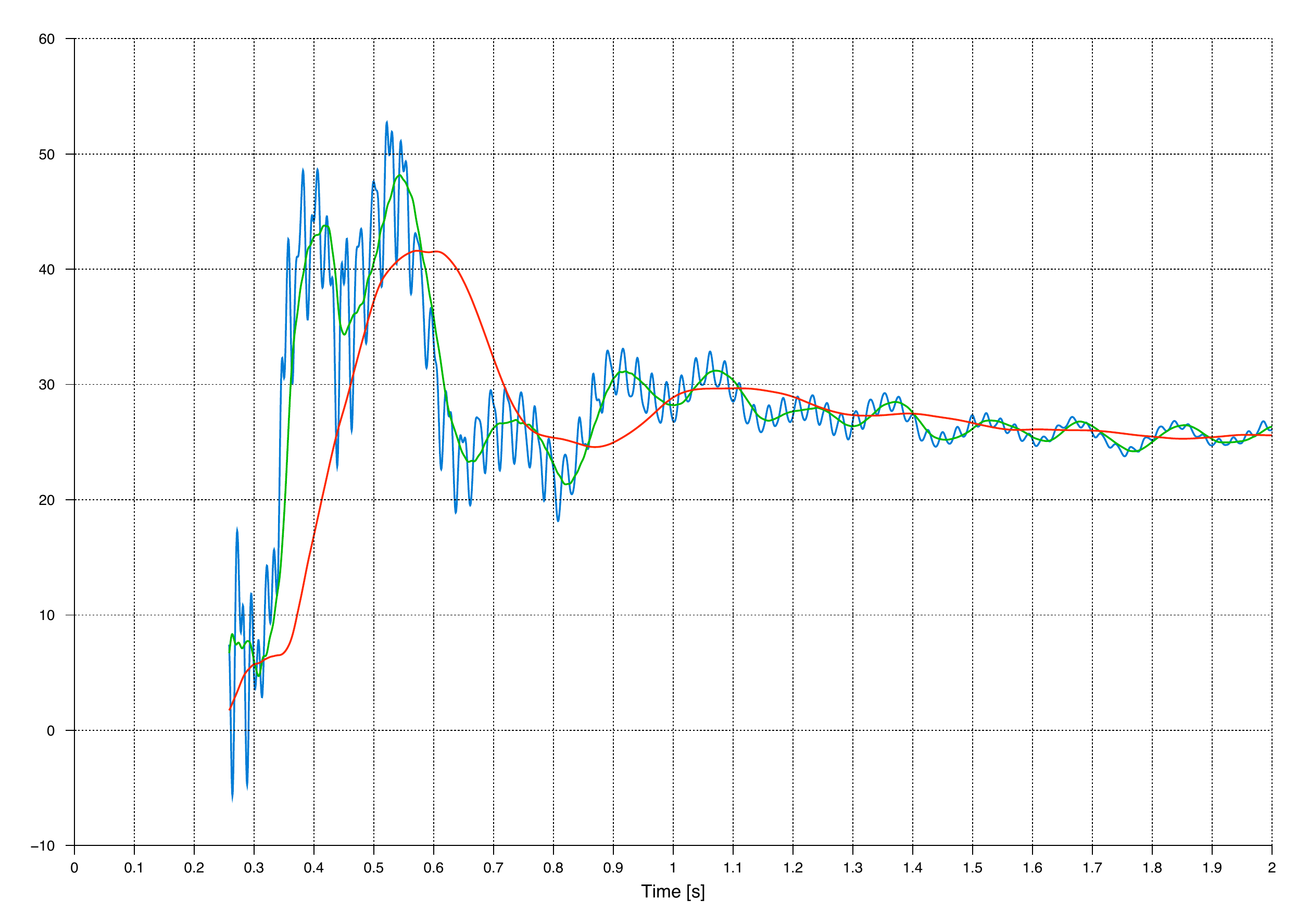}
\end{center}
\caption{Opening shock measurements filtered with moving averages
of different lenghts. The red curve obviously does not accurately
represent the forces during parachute opening.
\label{openingshockmovingaverages}}
\end{figure}
\section{How transients can spoil measurements}
In her Matura thesis \cite{carmen}, Carmen M\"uller tried to measure the opening shock
of small parachutes in the Ackeret wind tunnel of the Federal Institute of
Technology in Z\"urich (ETHZ). A specially constructed mounting tube
ejected the parachute using a piston powered by a small black powder charge
into the air stream.
While this accurately reflected the way parachutes are usually ejected
from small rockets, it also induced heavy oscillations in the measurement
equipment and probably also excited waves within the wind tunnel. As a
result, the data was almost unusable. The large amplitude transients
were at least one order of magnitude larger than the effects to be measured
(figure \ref{raw}).

Upon inspection of the data, two oscillatory components with frequencies
of about 40Hz and 12.2Hz were found. Assuming that the measured signal
is a superposition of the actual force to be measured $f(t)$ and a
transient of period $T$, convolution with a function of constant value
$T^{-1}$ on the interval $[-T,0]$, or equivalently a moving average over
the samples that fit in an interval of length $T$, eliminated the
transient. In this way, the two dominant transients could be removed,
and an approximation to the opening shock force was obtained.
The price to pay was of course that the signal $f$ was also
modified by the convolution. The modifications were large enough that it became
unreasonable to draw any quantitative conclusions from the measurements
(figure \ref{openingshockmovingaverages}).

Using an exponential weighting function, it was possible to slightly
improve the filtering, taking into account the damping of the transients.
Since the carrier of the convolution of the intervals of length
24ms and 83ms corresponding to the frequencies of the transients is
as an interval of at least 107ms, the convolution reduced the effective
resolution of the measurements to less than 10 samples per second, which
is of the same order of magnitude as the force changes during the
opening shock.

While this approach was successfull und fully adequate for a Matura
thesis, the question was still open whether a better filtering method
could reveal more details about the opening shock. Ideally, such a
filter would have a scaling parameter that allows to tune it
to different resolutions, where of course fine grained resolution
comes at the price of additional noise and incomplete filtering of
the observed transients. Below we show that a Kalman filter is ideally
suited for this situation, and that the system errors represent a set
of parameters that allow this kind of tuning of the filter.

\section{How to get rid of transients}
In the construction below, we will not develop the Kalman filter
equations, but only indicate the system equations from which the
filter can be derived. Our notation closely follows standard texts
on the discrete Kalman filter like \cite{catlin} or \cite{gelb}.
\subsection{The harmonic oscillator}
The measurement system in the wind tunnel works as follows. The
forces on the parachute and the ejection tube deflect the mounting
hardware by a tiny fraction measured by piezoelectric sensors.
So the mass $m$ of the mounting hardware is subject to two forces,
the forces we want to measure and the elastic forces, and we measure
its position $x$ as a function of time. If $D$ is Hook's constant for
the elastic forces, $x$ obeys the differential equation
\begin{equation}
m\ddot x+\alpha\dot x+Dx=f,
\label{deg2equation}
\end{equation}
where $f$ is the force we want to measure. In a pseudostatic situation,
we can assume that $\dot x= \ddot x =0$, and so the force is proportional
to the deflection value we measure. In the dynamic situation, we have
to infer $f(t)$ from $x(t)$ as a function of time.

A naive approach would use the measurements of $x$ and compute 
approximations of $\dot x$ and $\ddot x$ as finite difference quotients.
Then an approximation for $f$ can be computed from (\ref{deg2equation}).
However, numerical differentiation tends to be noisy, and the second
derivative will probably be almost useless. The problem gets even more serious
in light of the fact that in real applications, we need to filter several
times, so we need a method to estimate $f$ that does so with a controlled
level of noise.

This is the classical setting for a discrete Kalman filter. At any point in
time we would like the best estimate of the current position and speed of the
mass $m$ and the force $f$ based on measurements of $x$ at discrete points
$n\Delta t$ in time.
We use $x_n = (x(n\Delta t), \dot x(n\Delta t), f(n\Delta t))$ as the
state vector of a linear system. The time evolution is given by the 
differential equation (\ref{deg2equation}), which can be written,
using the familiar substitution $v(t)=\dot x(t)$  in
vector form as
\begin{equation}
\frac{d}{dt}\begin{pmatrix}x(t)\\ v(t)\end{pmatrix}
=
\begin{pmatrix}
0&1&0\\
-\frac{D}{m}&-\frac{\alpha}{m}&\frac{1}{m}
\end{pmatrix}
\begin{pmatrix}
x(t)\\v(t)\\f(t)
\end{pmatrix}
\end{equation}
Since we don't know anything about $f$, we assume for the purposes
of the Kalman filter, that $f$ is constant, i.~e.~$\dot f=0$.
Using the state vector $s=(x,v,f)$, the differential equation becomes
$$\frac{d}{dt}s=As$$
where
$$A=\begin{pmatrix}
0&1&0\\
-\frac{D}{m}&-\frac{\alpha}{m}&\frac{1}{m}\\
0&0&0
\end{pmatrix}
$$
Time evolution for a short time interval $\Delta t$ is given by
the matrix exponential
$$\varphi = e^{A\Delta t}.$$
The measurement equations are also very simple. Since we only measure
the position, we can use the measurement matrix
$$H=\begin{pmatrix}1&0&0\end{pmatrix}.$$

To completely define the Kalman filter, we have to specify the 
system and measurement errors.
The measurement error is a well known characteristic of the measurement
equipment. The system errors for $x$ and $v$ have to be kept very small,
as we believe the harmonic oscillator equation is followed quite
accurately. In contrast, the evolution equation for the force is only
a very rough approximation. We attribute almost all deviations of the real
system from the model to external forces. Thus the flucutations due to
measurement errors are the only sources of errors for $x$ and $v$,
while for $f$ we have to accept errors as large as the slope of $f$ times
the time step. If $\sigma_x^2$ is the position measurement variance and $\sigma_f^2$
the variance of the force,
then a reasonable approximation of the 
system error
covariance matrix $Q$ and the measurement error covariance $R$ is
$$Q=\begin{pmatrix}
\sigma_x^2&0&0\\
0&\sigma_x^2&0\\
0&0&\sigma_f^2\end{pmatrix},
\qquad
R=\begin{pmatrix}\sigma_x^2\end{pmatrix}.
$$
\begin{figure}
\begin{center}
\includegraphics[width=\hsize]{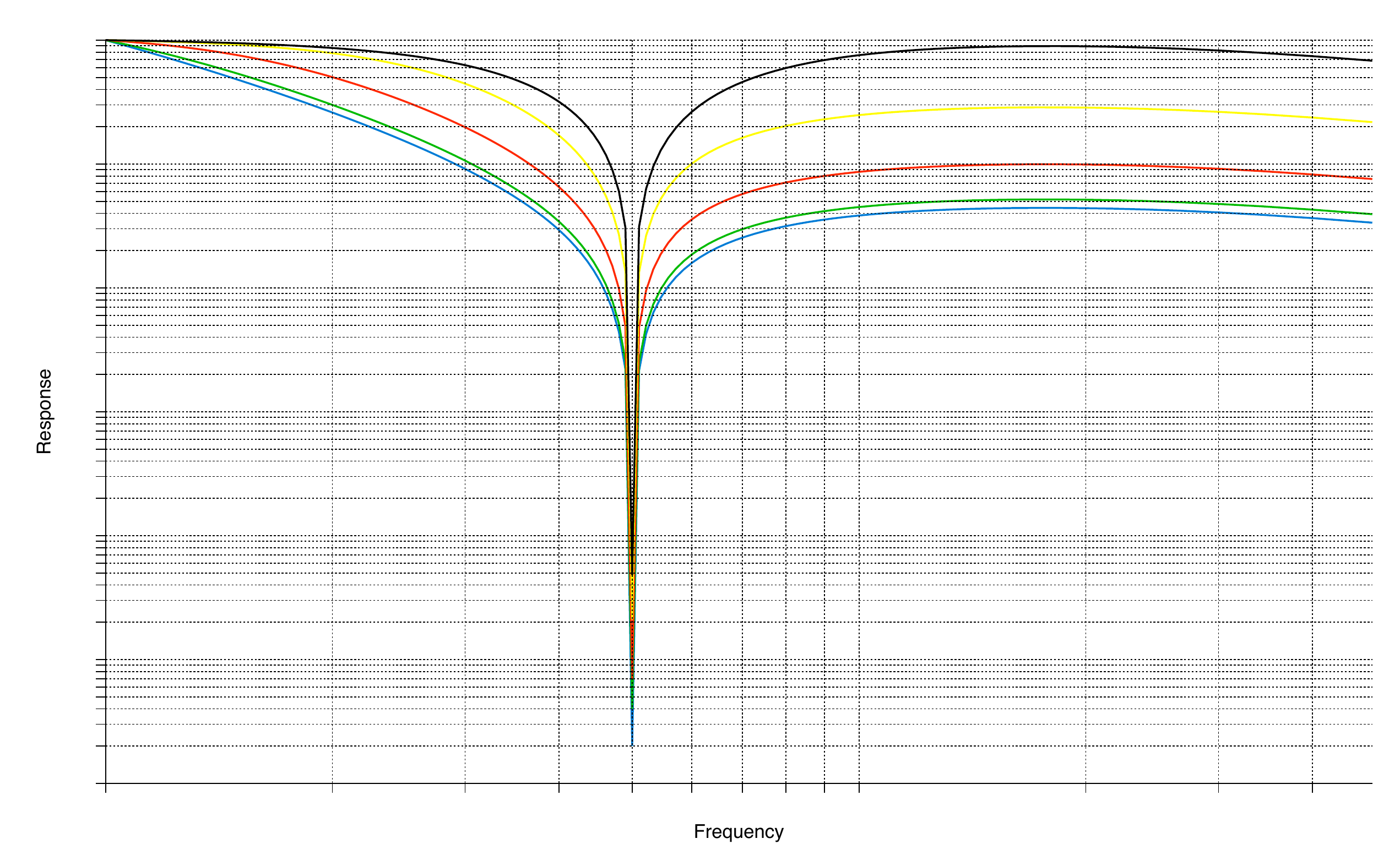}
\end{center}
\caption{Frequency response of Kalman filter based transient filter
for various values of $\sigma_f^2$\label{response}}
\end{figure}

In the example of the wind tunnel measurements of the opening shock,
the measurement error was about $0.001\text{N}^2$, while the
force error was on the order of $10^5\text{N}^2$.

Using a large value of $\sigma_f^2$ causes the Kalman filter to ``trust''
the measurements and agressively attribute deviations from the evolution
given by the differential equation to changing external forces. Of course,
$f$ will be much more noisy in this case, since measurement noise is converted
to a large degree to noise in $f$.
If $\sigma_f^2$
is small, then the filter tries not to modify $f$ too quickly, and
attributes deviations from the differential equation to a larger extent
to the measurement errors. The resulting $f$ will be much smoother,
but steep changes in $f$ will be filtered away.
Thus $\sigma_f^2$ is a natural scaling parameter we asked for in the
introduction.

\subsection{Frequency response}
Figure \ref{response} shows the numerically computed frequency
response for various values of $\sigma_f^2$.
As expected, the notch at 50Hz is very pronounced.
For large enough values of $\sigma_f^2$, there is hardly any loss for high
frequencies. 
For small values of $\sigma_f^2$, however, the loss is much larger,
rather similar to the $\sin x/x$ characteristic of the moving average
filter.

\subsection{Implementation considerations}
\subsubsection{Building the filter from $\omega$}
None of the coefficients $m$, $\alpha$ and $D$ of the differential
equation are known. Only the frequency and the damping of the transients
can be measured. The characteristic equation of the harmonic oscillator
equation
\begin{equation}
\omega^2+\frac{\alpha}m\omega+\frac{D}m=0
\label{charequation1}
\end{equation}
has two complex conjugate roots $\omega$ and $\bar\omega$.
To simplify the notation, we write $a=\frac{\alpha}D$ and $b=\frac{D}m$.
We conclude that
$$\operatorname{Re}\omega=\frac{a}2
\quad
\text{and}
\quad
\operatorname{Im}\omega=\sqrt{b-\frac{a^2}4}.$$
Given $\omega$, the differential equation 
can be obtained from
$$a=2\operatorname{Re}\omega
\quad
\text{and}
\quad
b=(\operatorname{Im}\omega)^2+\frac{a^2}4.$$
This allows to easily construct the Kalman filter for the oscillator
from frequencies obtained e.~b.~by spectral analysis.
\subsubsection{Precomputing the Kalman gain matrix}
It is well known that when all the matrices defining the Kalman filter
are constant, then the Kalman gain matrix converges. For given
values of the system and measurement errors, the Kalman gain
matrix can thus be computed beforehand and ``hard wired'' into the filter.

\subsection{Coupled oscillators}
If a system exhibits several oscillations, we can still assume that
the system is linear, and that each oscillation corresponds to a
harmonic oscillator with that frequency. Such a system can be described
by a higher dimensional state vector, and it is possible to construct
a Kalman filter for it.
However, it would become rather large, and the following simpler
approach is preferable. 

Each one of the coupled oscillators is accurately modelled by a
harmonic oscillator on which two external forces act, the force 
$f$ and the combined forces excerted by the other oscillators.
The latter have a well known frequency, so filtering the signal
in turn for all the frequencies of transients identified using
the method described above removes all the transients from the signal
$f$.

\section{Success stories}
The method described in the previous section was first applied to synthetic
data to verify that the filter is indeed capable of retrieving the synthetic
signal. Then the data was applied to the opening shock data that
originally motivated this research, and to a related problem involving
thrust curves of rocket motors.

\subsection{Test data}
\begin{figure}
\begin{center}
\includegraphics[width=0.8\hsize]{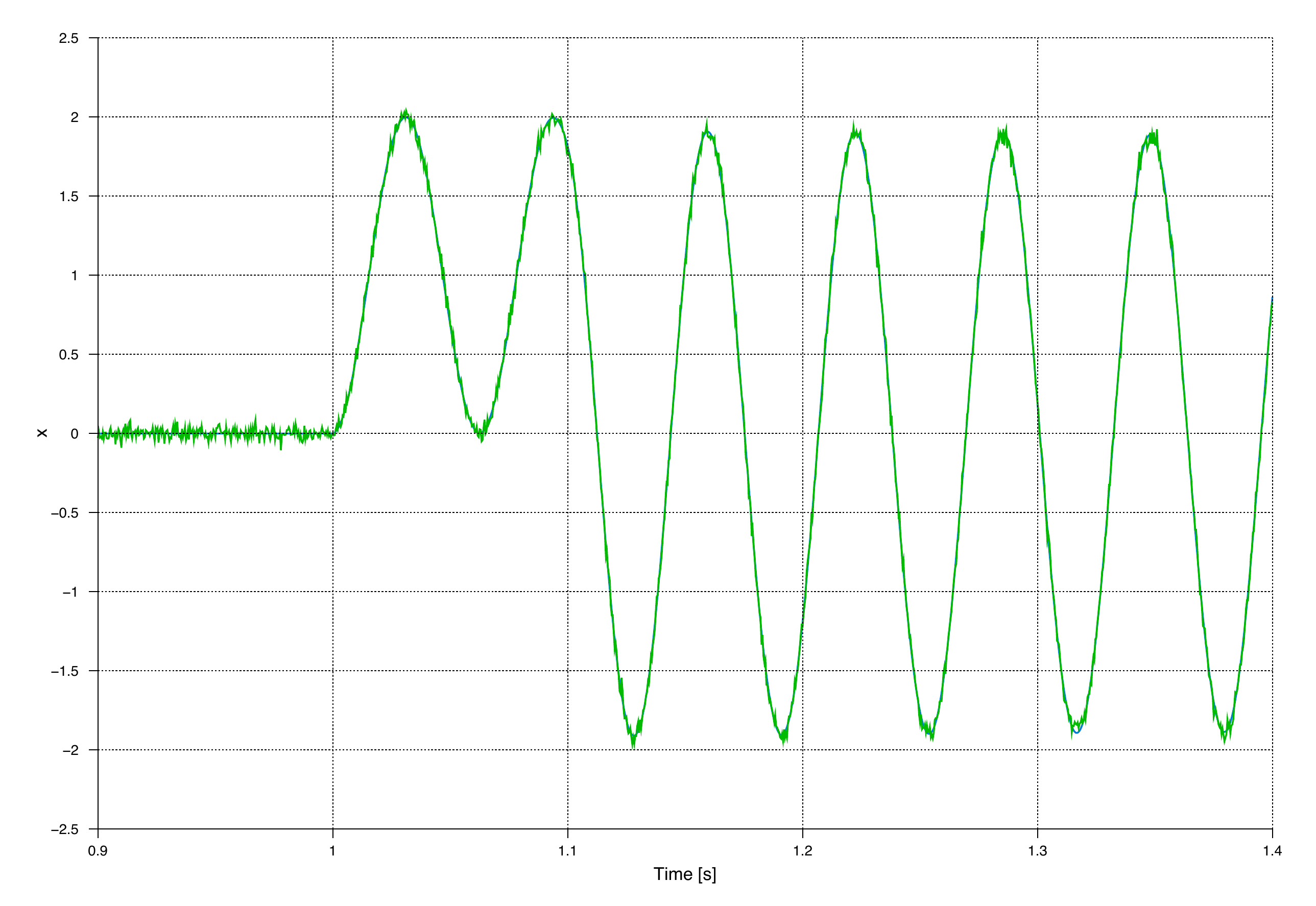}
\end{center}
\caption{Harmonic oscillator excited by a pulse of 0.1s duration at time $t=1$.\label{testdata}}
\end{figure}
\begin{figure}
\begin{center}
\includegraphics[width=0.8\hsize]{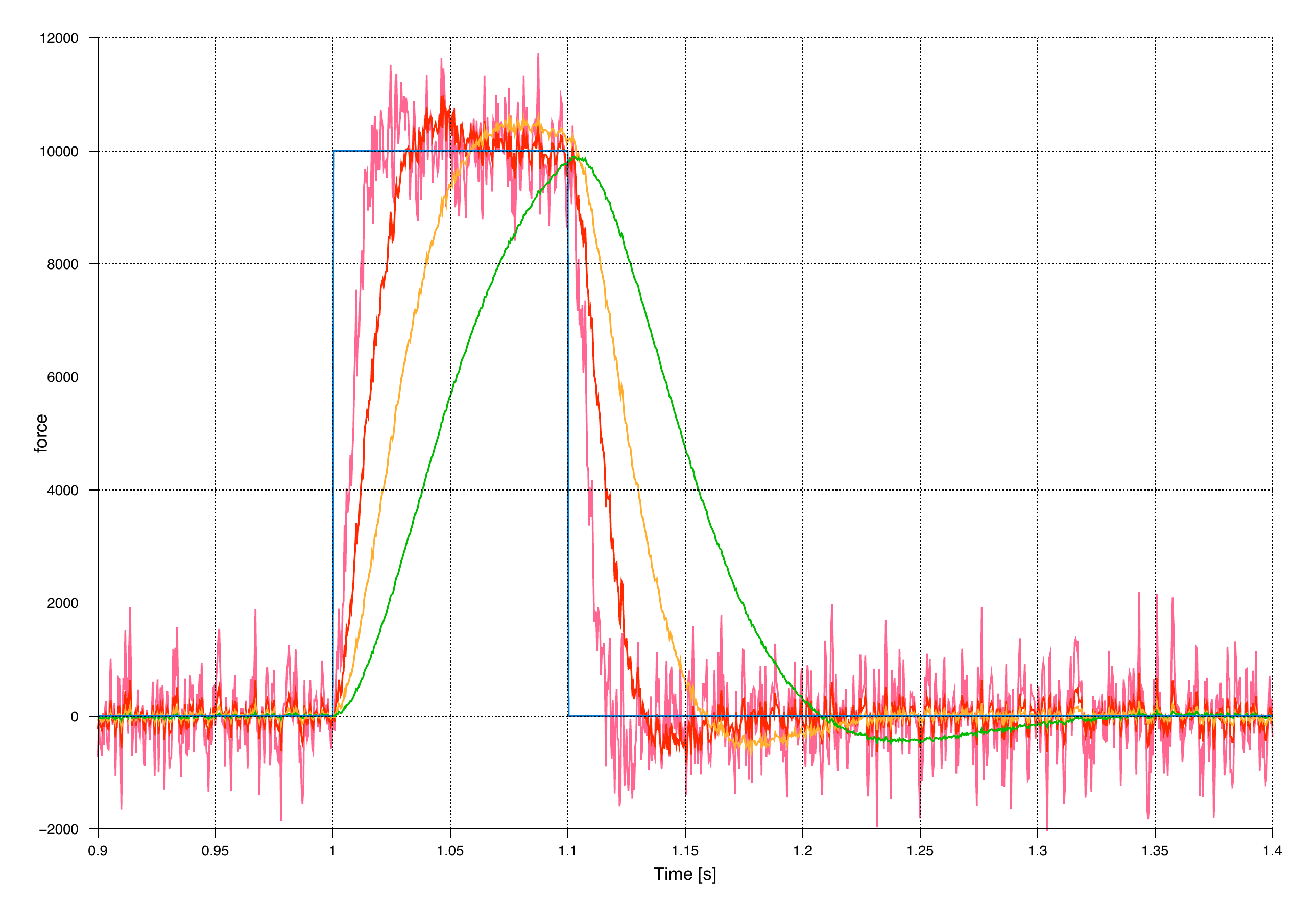}
\end{center}
\caption{Signal from figure \ref{testdata} filtered by a Kalman filter for increasing $\sigma_f^2$, compared to the original excitation (blue curve).\label{testfiltered}}
\end{figure}

To verify the properties of the Kalman filter based deconvolution of
transients, a solution for the harmonic oscillator equation
$$\ddot x + 0.1\cdot \dot x + 1000\cdot x = f(t)$$
was computed numerically for the force function
$$f(t)=\begin{cases}
0&\qquad t<1\vee t >1.1\\
10000&\qquad t\in[1,1.1]
\end{cases},
$$
and some gaussian noise with $\sigma=0.1$ added to simulate the measurement errors.
The resulting measurement data is shown in figure \ref{testdata}.
Figure \ref{testfiltered} shows the output of a Kalman filter based
deconvolution with $\sigma_f^2=10^{k}$, $k=3,4,5,6$.
As expected, noise increases with increasing $\sigma_f^2$. However, a simple
moving average over just a few samples can remove large parts of the noise
without distorting the signal very much.

\subsection{Opening shock of parachutes}
\begin{figure}
\begin{center}
\includegraphics[width=0.8\hsize]{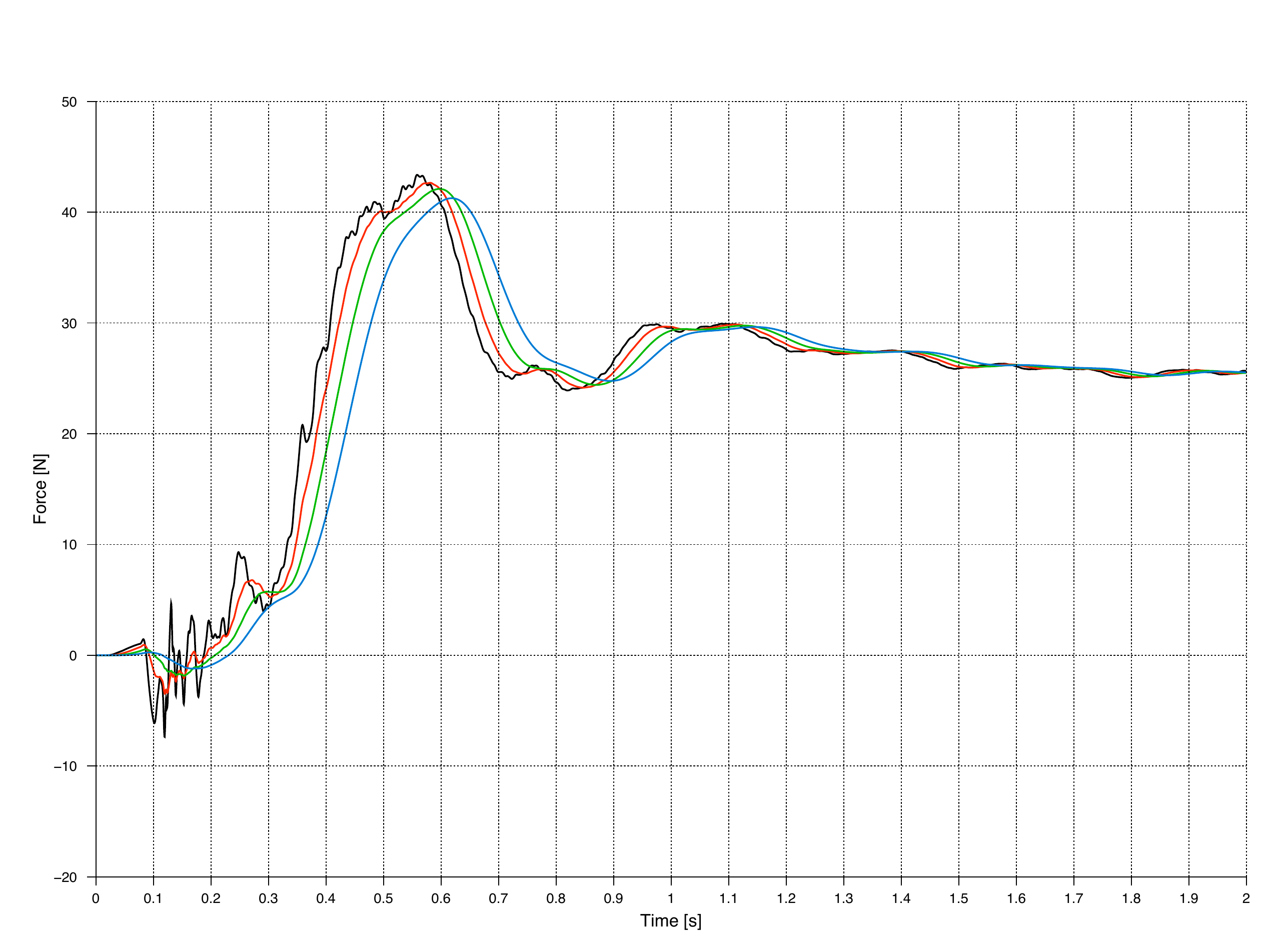}
\end{center}
\caption{Opening shock measurements from figure \ref{raw}
filtered with three Kalman filters for
different values of the force error variance $\sigma_f^2$.
\label{openingshockfiltered}}
\end{figure}
As a practical test case, the filtering method was applied to the opening
shock data. A chain of three kalman filters was constructed using the
following algorithm.
\begin{enumerate}
\item Take the raw as the current signal to filter.
\item Find the dominant frequency $\omega$ in the current signal.
\item Construct a Kalman filter for frequency $\omega$.
\item Filter the current signal, increase $\sigma_f^2$ until a new
dominant oscillation appears, in which case you continue at step 2,
or the signal becomes implausible.
\item Reduce $\sigma_f^2$ until you again get a physically plausible
signal.
\end{enumerate}
During this process, in addition to the frequencies already found in
\cite{carmen}, an additional oscillation at 85.7Hz was found. The
result of the filter for different values of $\sigma_f^2$ is shown
in figure \ref{openingshockfiltered}.
The curves following most closely the physical force signal were
obtained with $\sigma_f^2 = 100000 \text{N}^2$ and
$\sigma_f^2=200000\text{N}^2$, the second again shows some probably
unphysical transients.

The reconstruction of the signal using the Kalman filter at least
allows now to draw some quantitative conclusions from the opening
shock data.
E.~g.~it allowed to confirm the empirical rule used in some circles
that the opening shock is about
twice as strong as the drag of the fully inflated parachute.

\subsection{Thrust curve of high thrust rocket motors}
\begin{figure}
\begin{center}
\includegraphics[width=\hsize]{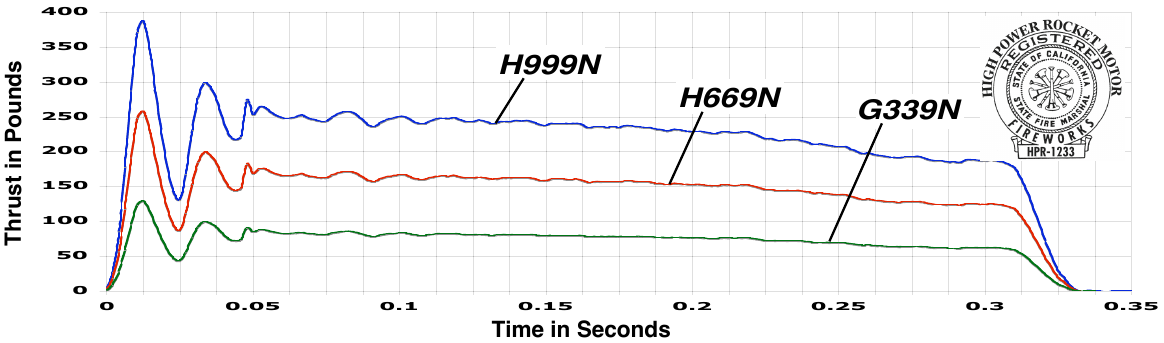}
\includegraphics[width=\hsize]{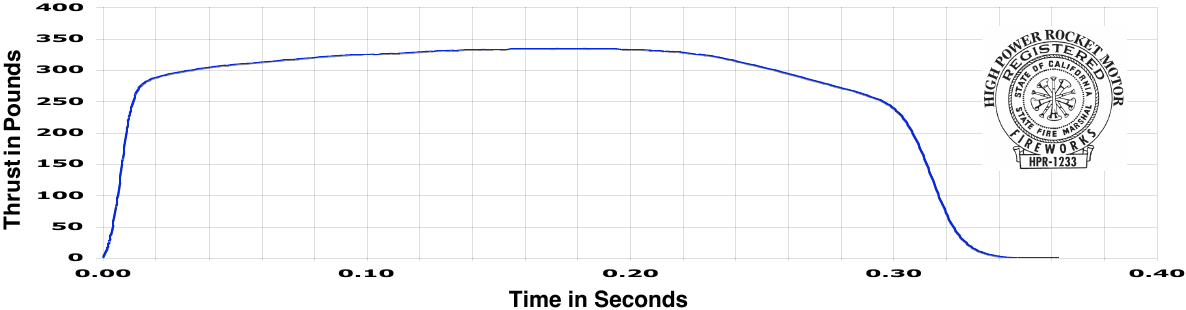}
\end{center}
\caption{Thrust curves of 38mm Warp-9 propellant rocket motors as they appear
on the instruction leaflet accompanying the reload kits. H999 top curve in top
diagram, I1299 in bottom diagram.\label{h999curve}}
\end{figure}
\begin{figure}
\begin{center}
\includegraphics[width=0.8\hsize]{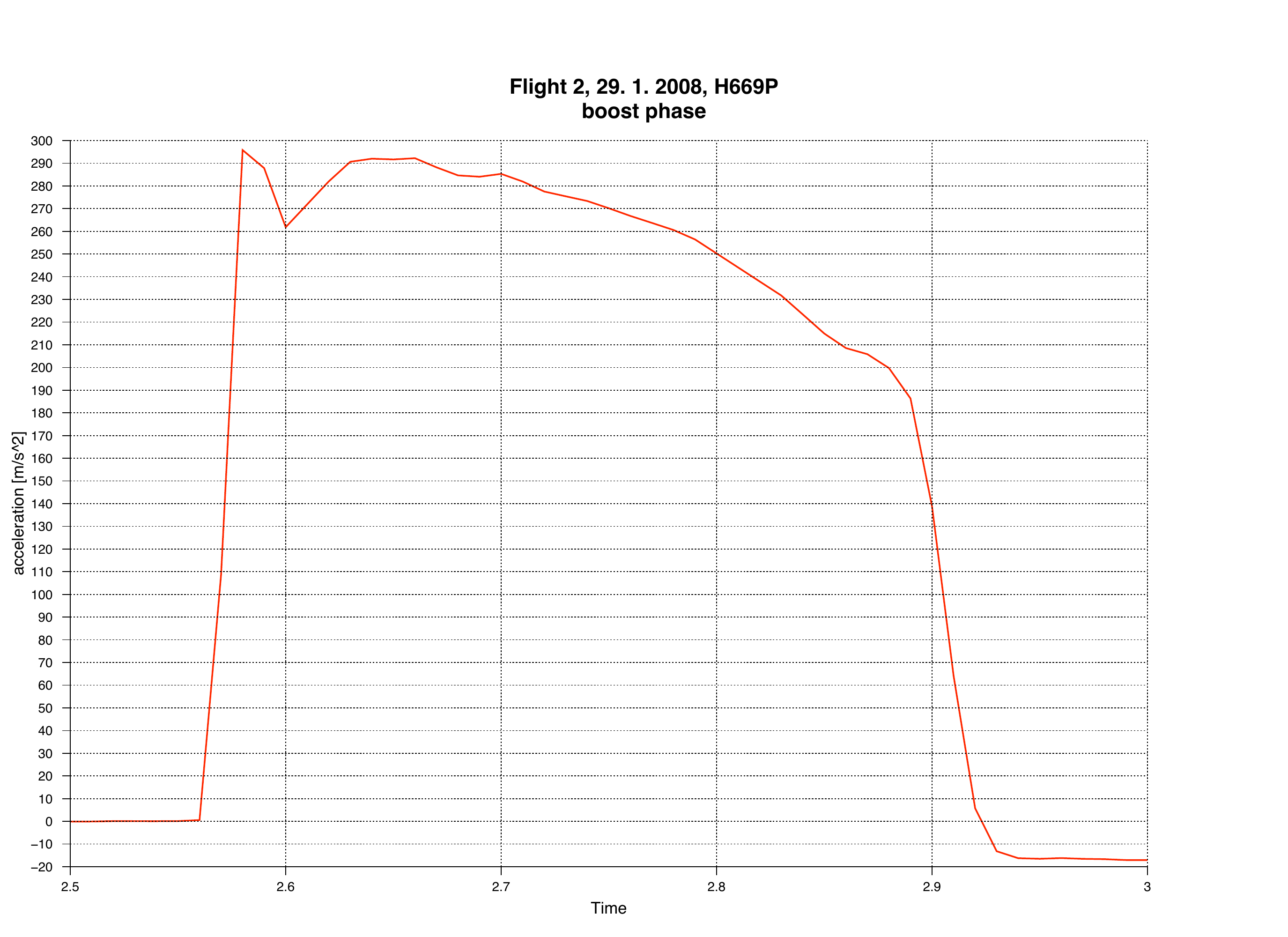}
\end{center}
\caption{Acceleration profile of a rocket flight on a H669 motor,
no oscillations\label{figurerocketware}}
\end{figure}
\begin{figure}
\begin{center}
\includegraphics[width=0.8\hsize]{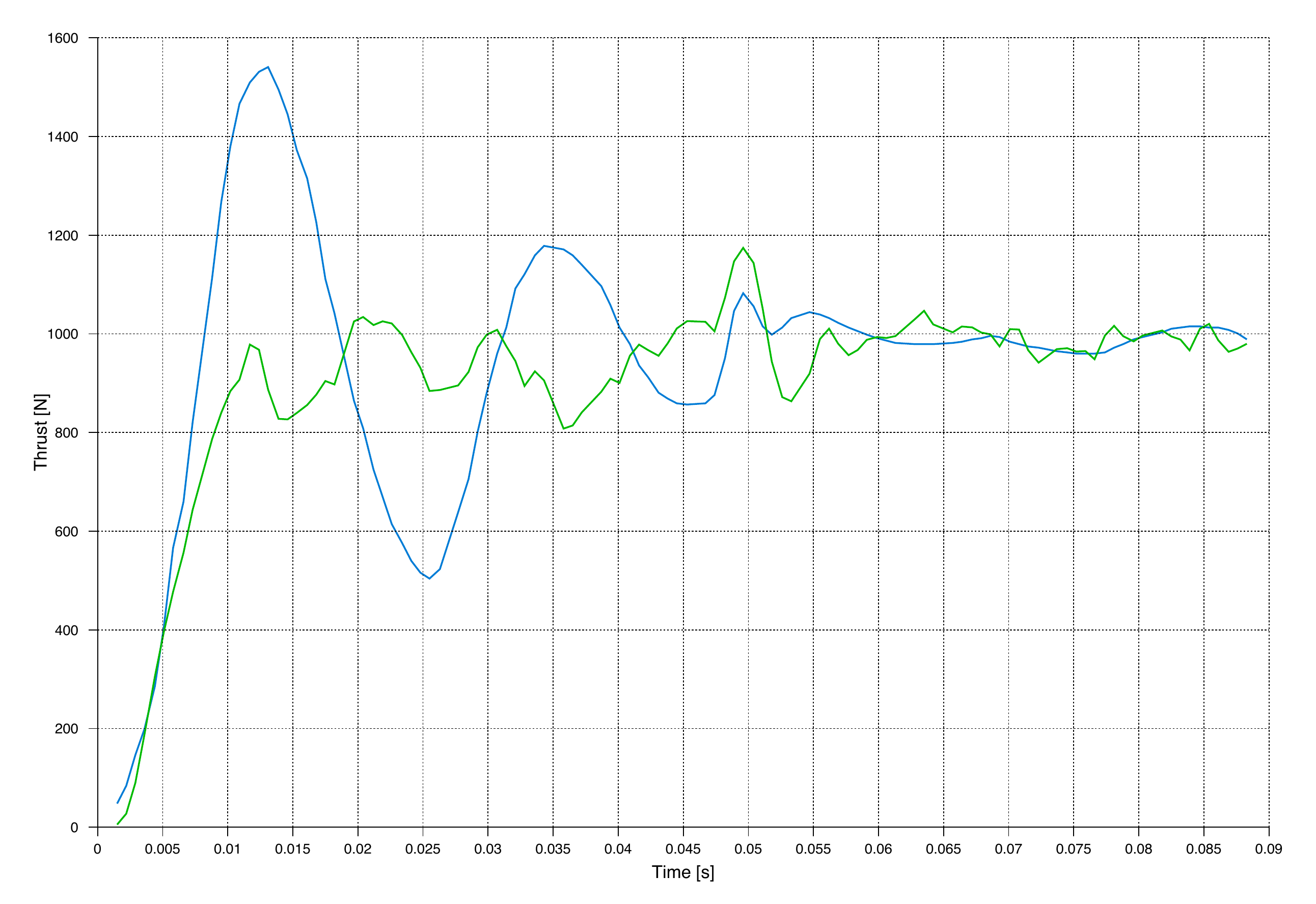}
\end{center}
\caption{Original thrust curve (blue) and version filtered for transients using a Kalman filter (green)\label{thrustfiltered}}
\end{figure}
Aerotech Consumer Aerospace produces rocket motors for hobby rockets
and small research rockets.
For each motor, thrust curves are available from various sources on
the Internet, and also from Aerotech.
The more recent motors using on the Warp-9 propellant have extremely
high thrust for a very short time interval.
The motor H999 has average
thrust of $999\text{N}$ during only $0.3\text{s}$.
The thrust curve shown in figure \ref{h999curve} bottom (from \cite{aerotechh999})
shows an initial
oscillation which cannot be explained by any known mechanism in
solid propellant motors. As a matter of fact, measuring the 
acceleration during an actual rocket flight with this motor, as was
done with the high precision altimeter developed during the rocketware
project (http://rocketware.othello.ch), shows no such oscillation
(see figure \ref{figurerocketware}). Also the only slightly larger
Motor I1299 (figure \ref{h999curve} bottom, from \cite{aerotechi1299})
has a completely flat thrust curve.
So most probably the oscillation is an artefact of the measuring equipment
used for these curves.
Of course, measuring thrust of a rocket motor is similar
to measuring aerodynamic forces in a wind tunnel. Thrust acts on the motor,
which is displaced slightly, the displacement is proportional to the force.
But without sufficient damping, the test stand holding the motor may oscillate
adding an oscillatory transient.

To corroborate this hypothesis,
the Kalman filter approach was applied to
a manually digitized version of the thrust curve. Figure~\ref{thrustfiltered}
shows the filtered version which is obviously much closer to the thrust
curve measured by the rocketware accelerometer.
The remaining noise was probably introduced when
manually reading thrust values from the printed thrust curve.

It turns out that the thrust curves were in fact measured on different
test stands. Thrust curves for I1299 were obtained on Aerotechs 
test stand, the test stand of the Tripoli Rocketry Association was
used for motors that show strong oscillations.

\section{Conclusion}
The results of actual measurement data filtered using the Kalman filter
approach show that oscillatory transients of at last the magnitude
of the signal can be effectively filtered to reveal the signal.
The measurements must be precise enough to render the full dynamic range,
otherwise the noise introduced by the filter will quickly become larger
than the signal. If this condition is satisfied, however, the method
can be expected to be close to optimal, as the Kalman filter has such
a property. If better models for $f$ are available, it is possible that
the method could be refined to yield even better results.

\bibliography{biblio}
\bibliographystyle{plain}

\end{document}